\documentclass[a4paper,10pt]{article}
\usepackage{amsmath}
\usepackage{amssymb}
\usepackage{graphics}
\usepackage{rotating}
\usepackage{cite}
\usepackage{color}

\DeclareMathOperator{\arctanh}{arctanh}

\title{ON SOME MODIFICATIONS OF THE FUETER OPERATOR}
\author{Daniel Alay\'{o}n-Solarz \underline{(danieldaniel@gmail.com)}}

\begin{document}

\maketitle

\begin{abstract}
We present some classes of functions that are defined on the quaternions as solutions for a
linear operator based on the Fueter operator. Unlike the Fueter regular functions;
in this case the identity function is a solution and solutions are closed under the quaternionic product. These classes  are a non-trivial extension of the Complex Holomorfic functions. In particular one of these classes are shown to satisfy a Cauchy-Riemann-like condition based on the spherical coordinates. Properties of these Classes as well as the generation of Fueter-regular functions are discussed.
\end{abstract}
\section*{I. INTRODUCTION}
"Quaternionic analysis" refers to many different types of extensions of the Cauchy-Riemann equations on the quaternion field. Fueter's approach; the so-called \textbf{regular functions}, is a rich theory that contains quaternionic versions of Cauchy's theorem and Cauchy's integral formula. Probably the best modern reference for the regular functions, in their general form, is Sudbery. The elegance of this generalization is perhaps clouded by the lack of some algebraic properties of the Complex Analysis: Regular functions cannot be multiplied or composed to obtain new regular functions and the identity function is not even regular. In recent years S.L Eriksson and H. Leutwiler have introduced a modified Dirac equation (related to the Cauchy-Riemann-Fueter condition for regularity), which gives rise to the Hyperholomorphic functions, which has also shown to be of great richness. It is interesting to note that Hyperholomorphic functions contains the identity function as well as its powers. Following that spirit this work is an attempt to study some other modifications of the Fueter operator whose solutions are algebraic closed (under the quaternionic product) and thus complex-like. What is complex-like in the quaternionic field is \textit{per se} a broad question. Another way to reformulate it could be the following: In how many ways the quaternions \textit{contains} the complex plane? The answer seems to depend on which special property of the complex numbers one is interested in. As we want a system of solutions for some Cauchy-Riemann equations that: are algebraic closed as quaternions and that contains the identity, we start by adding a rather strict condition: the functions evaluated at some point will not change the direction of the imaginary part of the quaternion. The study of such functions was first done by Rinehart and Cullen. Informally one obtains a process that turns the complex function $z^{n}$ onto the quaternion function $p^{n}$. We consider the set of commutative function such that they are analytic when restricted to some complex plane into the quaternions and it was introduced by S. De Leo and P. Rotelli. This is our starting point for the first of three classes of functions. The second class appears naturally when one evaluates the Fueter operator on the power of quaternions and appeared first on C.A Deavours article on the quaternion calculus. The third class, the most basic, are actually the Cullen functions. Deavours deduced that the Cullen functions would be solutions for the operator that defines the second class. We show that this result is incomplete as there exists solutions for this class that are essentially different from the Analytic Intrinsic functions studied by Cullen. Our main tool is the observation that the Fueter operator is invariant under some coordinate system related to the spherical coordinates. We show how our solutions classes are related to the Fueter's regularity. The bridge, in this case, is first based on a observation by Rinehart, that gives sufficient and necessary conditions for a complex function to be regular when turned into a quaternionic function by the Cullen method. These complex functions satisfy a non-standard version of the Cauchy-Riemann conditions.

%%%%%%%%%%%%%%%%%%%%%%%%%%%%%%%%%%%%%%%%%%%%%%%%%%%%%%%%%%%%%%%%%%%%%%%%%%%%%%%
%                                SECTION II
%%%%%%%%%%%%%%%%%%%%%%%%%%%%%%%%%%%%%%%%%%%%%%%%%%%%%%%%%%%%%%%%%%%%%%%%%%%%%%%

\section*{II. DEFINITIONS}
Let $p$ be a quaternion, we write $p$ in the canonical coordinates as:

\begin{equation}
p:=t+xi+yj+zk
\end{equation}

Let $f$ be a quaternionic map that commutes with its own argument in a open set $\omega$:

\begin{equation}
f(p)p = pf(p), \forall p \in \omega
\end{equation}

They are naturally related to complex maps. As we want to emphasize this we call these maps for \textit{Complex Extrinsic (CE)}.
\newtheorem{prop1}{Proposition}
\begin{prop1}
$f$ is EC (in some open set $\omega$) if and only if there exists real functions $u(p), v(p)$, defined on $\omega$ such that:
\begin{equation}
f(p) = u(p) + \iota v(p), \forall p \in \omega
\end{equation}
where
\begin{equation}
\iota := \frac{xi+yj+zk}{\sqrt{x^{2}+y^{2}+z^{2}}} \in S^{2}
\end{equation}
\end{prop1}

By fixing $\iota \in S^{2}$ we associate to $\omega$ a open subset $\tilde{\omega}$ in the upper complex plane,  in the following way:
\begin{equation}
z \in \tilde{\omega}_{\iota} \iff z = t + ri
\end{equation}
where
\begin{equation}
r:=\sqrt{x^{2}+y^{2}+z^{2}}
\end{equation}

In the same way we associate a CE map to a family of complex maps,
parametrized by $S^{2}$ and defined on $\tilde{\omega}_{\iota}$:
\begin{equation}
f_{\iota}(z) := u(z) + iv(z)
\end{equation}

In the rest of this work we will use the standard spherical coordinates:
\begin{equation}
\iota = (\cos\alpha \sin\beta, \sin\alpha \sin\beta, \cos \beta)
\end{equation}

We shall denote by $\iota_{\alpha}$, $\iota_{\beta}$  the derivative of $\iota$ respect to $\alpha$ and $\beta$ as a variable.

We say $f_{i_{p}}$ is a \textit{complex component} of the CE quaternionic map $f$. If $f$ is such that it only has one component we say $f$ is \textit{complex instrinsic (CI)}.

We add now differential conditions on our CE functions to relate them to complex analytic functions. All the Classes we are going to define are understood to be $C^{1}$.

We say a CE map $f$ is of \textit{Class I} if and only if it satisfies the the following Cauchy-Riemann equation:
\begin{equation}
\frac{\partial f}{\partial t} + \iota \frac{\partial f}{\partial r} =  0
\end{equation}
where
\begin{equation}
p = t + \iota r
\end{equation}

A CE function is called of \textit{Class II} if and only if:
\begin{equation}
\frac{\partial f }{\partial t}+ i\frac{\partial f}{\partial x} + j\frac{\partial f}{\partial y} + k\frac{\partial f }{\partial z} = \frac{-2v}{r}
\end{equation}

A function of Class I that is also CI, that is, with only one complex component is called of \textit{Class III}.

Finally, a CE function is called \textit{regular} if and only if:
\begin{equation}
\frac{\partial f }{\partial t}+ i\frac{\partial f}{\partial x} + j\frac{\partial f}{\partial y} + k\frac{\partial f }{\partial z} = 0
\end{equation}

The operator first defined in (11) is known as the Fueter (or sometimes Cauchy-Fueter or even  Cauchy-Riemann-Fueter) operator.  It has a left and right version:

\begin{equation}
\frac{\partial f}{\partial_{l}\bar{p}}=\frac{\partial f }{\partial t}+ i\frac{\partial f}{\partial x} + j\frac{\partial f}{\partial y} + k\frac{\partial f }{\partial z} 
\end{equation}

\begin{equation}
\frac{\partial f}{\partial_{r}\bar{p}}=\frac{\partial f }{\partial t}+ \frac{\partial f}{\partial x}i + \frac{\partial f}{\partial y}j + \frac{\partial f }{\partial z}k 
\end{equation}

We first rewrite $\frac{\partial f}{\partial_{l}\bar{p}}$ as:

\begin{equation}
\frac{\partial f}{\partial_{r}\bar{p}}=(\frac{\partial \bar{p}}{\partial t} )^{-1}\frac{\partial f}{\partial t}+ (\frac{\partial \bar{p}}{\partial x} )^{-1}\frac{\partial f}{\partial x}+ (\frac{\partial \bar{p}}{\partial y} )^{-1}\frac{\partial f}{\partial y} + (\frac{\partial \bar{p}}{\partial z} )^{-1}\frac{\partial f}{\partial z}
\end{equation}

Now let $\phi (t,x,y,z) = (t,r \cos\alpha \sin\beta,r  \sin\alpha \sin\beta, r \cos \beta)$ the parametrization. It is tempting to write then the following operator:

\begin{equation}
(\frac{\partial \bar{p}}{\partial t} )^{-1}\frac{\partial f((\phi (p))}{\partial t}+ (\frac{\partial \bar{p}}{\partial r} )^{-1}\frac{\partial f(\phi (p))}{\partial r }+ (\frac{\partial \bar{p}}{\partial \alpha} )^{-1}\frac{\partial f(\phi (p))}{\partial \alpha} + (\frac{\partial \bar{p}}{\partial \beta} )^{-1}\frac{\partial f(\phi (p))}{\partial \beta}
\end{equation}

And we inquire how this modified operator is related to the original Fueter operator. It is somehow surprising that they represent exactly the same operator:

\newtheorem{prop3}[prop1]{Proposition}
\begin{prop3}

\begin{equation}
 (\frac{\partial \bar{p}}{\partial t} )^{-1}\frac{\partial}{\partial t}+ (\frac{\partial \bar{p}}{\partial r} )^{-1}\frac{\partial}{\partial r }+ (\frac{\partial \bar{p}}{\partial \alpha} )^{-1}\frac{\partial}{\partial \alpha} + (\frac{\partial \bar{p}}{\partial \beta} )^{-1}\frac{\partial}{\partial \beta} =
\end{equation}
\begin{equation}
(\frac{\partial \bar{p}}{\partial t} )^{-1}\frac{\partial }{\partial t}+ (\frac{\partial \bar{p}}{\partial x} )^{-1}\frac{\partial }{\partial x}+ (\frac{\partial \bar{p}}{\partial y} )^{-1}\frac{\partial }{\partial y} + (\frac{\partial \bar{p}}{\partial z} )^{-1}\frac{\partial }{\partial z}  
\end{equation}
\end{prop3}

\textit{Proof} This is a mere use of the chain rule. We start writing (17) explicitly:

\begin{equation}
 (\frac{\partial \bar{p}}{\partial t} )^{-1}\frac{\partial}{\partial t}+ (\frac{\partial \bar{p}}{\partial r} )^{-1}\frac{\partial}{\partial r }+ (\frac{\partial \bar{p}}{\partial \alpha} )^{-1}\frac{\partial}{\partial \alpha} + (\frac{\partial \bar{p}}{\partial \beta} )^{-1}\frac{\partial}{\partial \beta} =
\end{equation}

\begin{equation}
 \frac{\partial}{\partial t}+ \iota \frac{\partial}{\partial r }-  r^{-1} \iota_{\alpha}^{-1} \frac{\partial}{\partial \alpha} -  r^{-1} \iota_{\beta}^{-1} \frac{\partial}{\partial \beta} =
\end{equation}

\begin{equation*}
\frac{\partial}{\partial t}+ \iota (\frac{\partial x}{\partial r }\frac{\partial }{\partial x }+\frac{\partial y}{\partial r }\frac{\partial }{\partial y} +\frac{\partial z}{\partial r }\frac{\partial }{\partial z }) - r^{-1} \iota_{\alpha}^{-1}(\frac{\partial x}{\partial \alpha }\frac{\partial }{\partial x }+\frac{\partial y}{\partial \alpha }\frac{\partial }{\partial y} +\frac{\partial z}{\partial \alpha }\frac{\partial }{\partial z })-
\end{equation*}

\begin{equation}
 r^{-1} \iota_{\beta}^{-1}(\frac{\partial x}{\partial \beta }\frac{\partial }{\partial x }+\frac{\partial y}{\partial \beta }\frac{\partial }{\partial y} +\frac{\partial z}{\partial \beta }\frac{\partial }{\partial z })  =
\end{equation}

\begin{equation*}
\frac{\partial}{\partial t}+(\iota \frac{\partial x}{\partial r } -  r^{-1} \iota_{\alpha}^{-1}\frac{\partial x}{\partial \alpha } -  r^{-1} \iota_{\beta}^{-1}\frac{\partial x}{\partial \beta })\frac{\partial }{\partial x} +(\iota \frac{\partial y}{\partial r } -  r^{-1} \iota_{\alpha}^{-1}\frac{\partial y}{\partial \alpha } -  r^{-1} \iota_{\beta}^{-1}\frac{\partial y}{\partial \beta })\frac{\partial }{\partial y}
\end{equation*}

\begin{equation}
+ (\iota \frac{\partial z}{\partial r } -  r^{-1} \iota_{\alpha}^{-1}\frac{\partial z}{\partial \alpha } -  r^{-1} \iota_{\beta}^{-1}\frac{\partial z}{\partial \beta })\frac{\partial }{\partial z}
\end{equation}

Finally, the validity of the proposition depends of the fact that:

\begin{equation}
\iota \frac{\partial x}{\partial r } -  r^{-1} \iota_{\alpha}^{-1}\frac{\partial x}{\partial \alpha } -  r^{-1} \iota_{\beta}^{-1}\frac{\partial x}{\partial \beta } = i
\end{equation}
\begin{equation}
\iota \frac{\partial y}{\partial r } -  r^{-1} \iota_{\alpha}^{-1}\frac{\partial y}{\partial \alpha } -  r^{-1} \iota_{\beta}^{-1}\frac{\partial y}{\partial \beta } = j
\end{equation}
\begin{equation}
\iota \frac{\partial z}{\partial r } -  r^{-1} \iota_{\alpha}^{-1}\frac{\partial z}{\partial \alpha } -  r^{-1} \iota_{\beta}^{-1}\frac{\partial z}{\partial \beta } = k
\end{equation}

%%%%%%%%%%%%%%%%%%%%%%%%%%%%%%%%%%%%%%%%%%%%%%%%%%%%%%%%%%%%%%%%%%%%%%%%%%%%%%%
%                                SECTION III
%%%%%%%%%%%%%%%%%%%%%%%%%%%%%%%%%%%%%%%%%%%%%%%%%%%%%%%%%%%%%%%%%%%%%%%%%%%%%%%
\section*{III. SOME PROPERTIES ON CLASSES I TO III}
Given two regular functions $f,g$ (not necesarilly CE) in some open set $\omega$ in general it is not true that $fg$ or $gf$ will be regular. This is an important difference between the Fueter regularity and the classical analyticity, and makes the generation of regular functions more difficult. Note that that in the Fueter regularity the identity function $f(p) = p$ is not regular. Classes I to III, on the other side, behave nicely with respect to the quaternionic product and contain the identity function.
\newtheorem{prop20}[prop1]{Proposition}
\begin{prop20}
Let $f$,$g$ be both of Class I,II or III then $fg=gf$ and $f+g$ are of respective Class I,II or III. If there
exists  the algebraic inverse $f^{-1}$ then it is also of respective Class I, II or III.
\end{prop20}

\textbf{Proof} This is a straight-forward calculation and will be omitted.

\newtheorem{prop4}[prop1]{Proposition}
\begin{prop4}
Class I $\supset$ Class II $\supset$ Class III
\end{prop4}
\textit{Proof} We first observe that the Class III is by definition contained in the Class I. Let $f=u(p)+\iota v(p)$ be of Class II, if we write $f = u + \iota v$ this is equivalent as to say that $u$ and $w:= \frac{v}{r}$ satisfy the following set of equations in $(t,x,y,z)$ coordinates:

\begin{eqnarray}
\frac{\partial u}{\partial t} - \frac{\partial w}{\partial x}x - \frac{\partial w}{\partial y}y - \frac{\partial v}{\partial z}z &=& w
\\
\frac{\partial u}{\partial x} + \frac{\partial w}{\partial t}x &=& \frac{\partial w}{\partial z}y - \frac{\partial w}{\partial y}z
\\
\frac{\partial u}{\partial y} + \frac{\partial w}{\partial t}y &=& \frac{\partial w}{\partial x}z - \frac{\partial w}{\partial z}x
\\
\frac{\partial u}{\partial z} + \frac{\partial w}{\partial t}z &=& \frac{\partial w}{\partial y}x - \frac{\partial w}{\partial x}y
\end{eqnarray}
Multiplying (27) by $x$, (28) by $y$ and (29) by $z$ and summing turns these three equations in only one. So we obtain only two equations, namely:
\begin{eqnarray}
\frac{\partial u}{\partial t} - \frac{\partial w}{\partial x}x - \frac{\partial w}{\partial y}y - \frac{\partial v}{\partial z}z &=& w \nonumber
\\
\frac{\partial u}{\partial x}x + \frac{\partial u}{\partial y}y + \frac{\partial u}{\partial z}z + \frac{\partial w}{\partial t}x^{2} + \frac{\partial w}{\partial t}y^{2} + \frac{\partial w}{\partial t}z^{2}  &=& 0
\end{eqnarray}
Turning back these two equations in coordinates $t,r,\iota$ is exactly (9). So $f$ is of Class I. Now let $f$ be of Class III.
Let $(\alpha, \beta)$ be the parametrization of $S^{2}$, in spherical coordinates. We observe that if $f$ is of Class III this then:
\begin{equation}
\frac{\partial u}{\partial \alpha} = \frac{\partial u}{\partial \beta} =
\frac{\partial v}{\partial \alpha} = \frac{\partial v}{\partial \beta} =
0
\end{equation}
Applying the Fueter operator to $f$ in $(t,r,\iota) = (t,r,\alpha, \beta)$ coordinates is:
\begin{equation}
\frac{\partial f}{\partial t} + \iota \frac{\partial f}{\partial r} + (\frac{\partial \bar{p}}{\partial \alpha})^{-1} \frac{\partial f}{\partial \alpha} + (\frac{\partial \bar{p}}{\partial \beta})^{-1} \frac{\partial f}{\partial \beta}
\end{equation}
We observe that
\begin{equation}(\frac{\partial \bar{p}}{\partial \alpha})^{-1} \frac{\partial f}{\partial \alpha} = - (\frac{\partial \iota}{\partial \alpha}r)^{-1}(\frac{\partial f}{\partial \alpha}) = (\frac{\partial \bar{p}}{\partial \beta})^{-1} \frac{\partial f}{\partial \beta} = - (\frac{\partial \iota}{\partial \beta}r)^{-1}(\frac{\partial f}{\partial \beta}) = -\frac{v}{r}
\end{equation}
A Class III functions is by definition of Class I, the first two summands in (32) are precisely the Class I condition and so they vanish, this together with the previous calculation shows that:
\begin{equation}
\frac{\partial f}{\partial t} + \iota \frac{\partial f}{\partial r} + (\frac{\partial \bar{p}}{\partial \alpha})^{-1} \frac{\partial f}{\partial \alpha} + (\frac{\partial \bar{p}}{\partial \beta})^{-1} \frac{\partial f}{\partial \beta} = \frac{-2v}{r}
\end{equation}
and so $f$ is of Class II.

By a straight-forward calculation is obtained a simple formula for the jacobian determinant of a Class I to III function at a given point.

\newtheorem{prop5}[prop1]{Proposition}
\begin{prop5}
Let $f$ be a function of Class II, the determinant of the jacobian is given by the following formula:
\begin{equation}
det(D_{p}f) = \vert\frac{\partial f(p)}{\partial t}\vert^{2} \frac{v(p)^2}{r^2}
\end{equation}
where:
\begin{equation}
\vert\frac{\partial f}{\partial t}\vert^{2} := (\frac{\partial u(p)}{\partial t})^{2} + (\frac{\partial v(p)}{\partial t})^{2}
\end{equation}
\end{prop5}

%%%%%%%%%%%%%%%%%%%%%%%%%%%%%%%%%%%%%%%%%%%%%%%%%%%%%%%%%%%%%%%%%%%%%%%%%%%%%%%
%                                SECTION IV
%%%%%%%%%%%%%%%%%%%%%%%%%%%%%%%%%%%%%%%%%%%%%%%%%%%%%%%%%%%%%%%%%%%%%%%%%%%%%%%

\section*{IV. CLASS II IS NOT CLASS III}
Let $f$ be a Class III function defined on some open set $\omega$. By definition it has only one complex component and this component is analytical. We denote by $\tilde{f}$ this analytical function. Conversely, given any analytic function $\tilde{f}$ defined in (some open set in) the upper complex plane and following Cullen we associate $f$ as the Class III function that has $\tilde{f}$ as its complex component. Deavours showed, by a different manner than we have used, that Cullen functions, that is, Class III functions will be of Class II. He stated that only functions generated in the above (Cullen) manner are of Class II. This would mean that the Class II and III are actually the same. We show that there exists an abundance of Class II solutions, at least defined locally, that are not Class III.

\newtheorem{prop6}[prop1]{Proposition}

\begin{prop6}
Let $f$ be a Class II function, defined in some open set, then $f$ for the spherical parametrization $(\alpha, \beta)$ of $S^{2}$ satisfy the following equations:
\begin{equation}
\frac{\partial v}{\partial \alpha }(\sin\beta)^{-1} + \frac{\partial u}{\partial \beta} = \frac{\partial u}{\partial \alpha }(\sin\beta)^{-1} - \frac{\partial v}{\partial \beta} = 0
\end{equation}
\end{prop6}
\textit{Proof} Because $f$ is assumed of Class II then for the parametrization of the sphere $f$ satisfies:
\begin{equation}
(\frac{\partial \bar{p}}{\partial \alpha })^{-1}\frac{\partial f}{\partial \alpha } + (\frac{\partial \bar{p}}{\partial \beta })^{-1}\frac{\partial f}{\partial \beta } = \frac{-2v}{r}
\end{equation}

We observe that

\begin{equation}
(\frac{\partial \bar{p}}{\partial \alpha})^{-1} = (\frac{\partial}{\partial \alpha }(t-r \iota))^{-1} =
(-r)^{-1}(\frac{\partial \iota}{\partial \alpha})^{-1}
\end{equation}
and
\begin{equation}
(\frac{\partial \bar{p}}{\partial \beta})^{-1} = (\frac{\partial}{\partial \beta }(t-r \iota))^{-1} =
(-r)^{-1}(\frac{\partial \iota}{\partial \beta})^{-1}
\end{equation}

By the assumption: $r \neq 0$ and we rewrite (38) as:
\begin{equation}
(\frac{\partial\iota}{\partial \alpha})^{-1}\frac{\partial f}{\partial \alpha } + (\frac{\partial\iota}{\partial \beta})^{-1}\frac{\partial f}{\partial\beta} = 2v
\end{equation}
We substitute $f=u+\iota v$ in (41) and obtain:
\begin{equation}
(\frac{\partial\iota}{\partial \alpha})^{-1}(\frac{\partial u}{\partial \alpha } + \iota\frac{\partial v}{\partial \alpha } + \frac{\partial \iota}{\partial \alpha }v) +
(\frac{\partial\iota}{\partial \beta})^{-1}(\frac{\partial u}{\partial \beta } + \iota\frac{\partial v}{\partial \beta } + \frac{\partial \iota}{\partial \beta }v) = 2v
\end{equation}
then
\begin{equation}
(\frac{\partial\iota}{\partial \alpha})^{-1}\frac{\partial u}{\partial \alpha } + (\frac{\partial\iota}{\partial \alpha})^{-1}\iota\frac{\partial v}{\partial \alpha } +
(\frac{\partial\iota}{\partial \beta})^{-1}\frac{\partial u}{\partial \beta } + (\frac{\partial\iota}{\partial \beta})^{-1}\iota\frac{\partial v}{\partial \beta } = 0
\end{equation}
Now we write (43) explicitly in terms of the parametrization:
\begin{equation}
\frac{\partial \iota}{\partial \alpha} = \left(\begin{array}{c} -\sin\alpha \sin\beta \\ \cos\alpha \sin\beta \\ 0 \end{array} \right)
\end{equation}
and
\begin{equation}
\frac{\partial \iota}{\partial \beta} = \left(\begin{array}{c} \cos\alpha \cos\beta \\ \sin\alpha \cos\beta \\ -\sin\beta \end{array} \right)
\end{equation}

the square norm of these vectors, are:
\begin{equation}
{\Vert\frac{\partial \iota}{\partial \alpha}\Vert}^{2}=
\sin^{2}\beta.
\end{equation}
\begin{equation}
{\Vert\frac{\partial \iota}{\partial \beta}\Vert}^{2} = 1
\end{equation}
Since $\iota^{2}= -1$ then $\frac{\partial \iota}{\partial \alpha}\iota + \iota \frac{\partial \iota}{\partial \alpha}=0$, and the same holds if we replace $\alpha$ by $\beta$. As these two imaginary vectors anticommute their quaternionic product is exactly their cross product. The anticommutativity still holds if we replace any of these vectors by their inverses.
We continue this calculation expressing these vectors in terms of the parametrization:

\begin{eqnarray}
(\frac{\partial \iota}{\partial \alpha})^{-1} = (\sin\beta)^{-1} \left(\begin{array}{c} \sin\alpha \\ -\cos\alpha \\ 0 \end{array} \right)
\\
(\frac{\partial \iota}{\partial \beta})^{-1} = \left(\begin{array}{c} -\cos\alpha \cos\beta \\ -\sin\alpha \cos\beta \\ \sin\beta \end{array} \right)
\\
(\frac{\partial \iota}{\partial \alpha})^{-1} \iota = (\sin\beta)^{-1} \left(\begin{array}{c} -\cos\alpha \cos\beta \\ -\sin\alpha \cos\beta \\ \sin\beta \end{array} \right)
\\
(\frac{\partial \iota}{\partial \beta})^{-1} \iota = \left(\begin{array}{c} -\sin\alpha \\ \cos\alpha \\ 0 \end{array} \right)
\end{eqnarray}
We substitute this expressions in (43):
\begin{equation}
(\sin\beta)^{-1}\left(\begin{array}{c} \sin\alpha \\ -\cos\alpha \\ 0 \end{array} \right)\frac{\partial u}{\partial \alpha}+
\left(\begin{array}{c} -\cos\alpha \cos\beta \\ -\sin\alpha \cos\beta \\ \sin\beta \end{array} \right)\frac{\partial v}{\partial \alpha}(\sin\beta)^{-1}+
\end{equation}
\begin{equation} \nonumber
+\left(\begin{array}{c} -\cos\alpha \cos\beta \\ -\sin\alpha \cos\beta \\ \sin\beta \end{array} \right)\frac{\partial u}{\partial \beta}+
\left(\begin{array}{c} -\sin\alpha \\ \cos\alpha \\ 0 \end{array} \right)\frac{\partial v}{\partial \beta} = 0
\end{equation}
Which finally implies:
\begin{eqnarray}
\frac{\partial v}{\partial \alpha}(\sin\beta)^{-1}+\frac{\partial u}{\partial \beta} = 0
\\
\frac{\partial u}{\partial \alpha}(\sin\beta)^{-1}-\frac{\partial v}{\partial \beta}= 0
\end{eqnarray}

Consider the function $f = xr^{-1}\iota$ is a Class I function that is not Class II.

We observe also that the following function:

\begin{equation}
\rho(\alpha,\beta):= \alpha + \iota \ln (\tan(\beta/2))
\end{equation}

is a Class II function that is not of Class III. The expression for this solution in $(t,x,y,z)$ coordinates is:

\begin{equation}
\rho(x,y,z) := \arctan( \frac{x}{y}) + \iota \arctanh(\frac{z}{\sqrt{x^{2}+y^{2}+z^{2}}})
\end{equation}

From which we construct two more solutions:

\begin{equation}
\varrho(x,y,z) := \arctan( \frac{y}{z}) + \iota \arctanh(\frac{x}{\sqrt{x^{2}+y^{2}+z^{2}}})
\end{equation}
\begin{equation}
\sigma(x,y,z) := \arctan( \frac{z}{x}) + \iota \arctanh(\frac{y}{\sqrt{x^{2}+y^{2}+z^{2}}})
\end{equation}

As expected, $\bar{\rho}$ denotes the quaternionic conjugate of the function $\rho$. We appreciate chirality by the following observation, which holds also for $\varrho$,$\sigma$:
\begin{equation}
\frac{\partial \rho}{\partial _{l} \bar{p}} = \frac{\partial \bar{\rho}}{\partial _{r} \bar{p}}
\end{equation}

This observation can be generalized:

\newtheorem{prop17}[prop1]{Proposition}

\begin{prop17}
Let $f$ be a left-Class II function such that

\begin{equation}
\frac{\partial u}{\partial t} = \frac{\partial u}{\partial r} = \frac{\partial v}{\partial t} = \frac{\partial v}{\partial r} = 0
\end{equation}

then the conjugate of $f$ will be a right-Class II.

\end{prop17}

\textit{Proof} We consider 
\begin{equation}
(\frac{\partial\iota}{\partial \alpha})^{-1}(\frac{\partial u}{\partial \alpha } + \iota\frac{\partial v}{\partial \alpha } + \frac{\partial \iota}{\partial \alpha }v) +
(\frac{\partial\iota}{\partial \beta})^{-1}(\frac{\partial u}{\partial \beta } + \iota\frac{\partial v}{\partial \beta } + \frac{\partial \iota}{\partial \beta }v) = 2v
\end{equation}

which can be rewritten as:

\begin{equation}
(\frac{\partial u}{\partial \alpha } + \iota\frac{\partial (-v)}{\partial \alpha } + \frac{\partial \iota}{\partial \alpha }(-v)) (\frac{\partial\iota}{\partial \alpha})^{-1}+
(\frac{\partial u}{\partial \beta } + \iota\frac{\partial (-v)}{\partial \beta } + \frac{\partial \iota}{\partial \beta }(-v))(\frac{\partial\iota}{\partial \beta})^{-1} = -2v
\end{equation}

On the other side, a Class III function is central:

\newtheorem{prop18}[prop1]{Proposition}
\begin{prop18}
Let $f$ be a Class II function. Then

\begin{equation}
\frac{\partial f}{\partial_{l} \bar{p}} =\frac{\partial f}{\partial_{r} \bar{p}} = \frac{-2v}{r}
\end{equation}

if and only if $f$ is of Class III
\end{prop18}
\textit{Proof} Suppose 

\begin{equation}
(\frac{\partial\iota}{\partial \alpha})^{-1}(\frac{\partial u}{\partial \alpha } + \iota\frac{\partial v}{\partial \alpha } + \frac{\partial \iota}{\partial \alpha }v) +
(\frac{\partial\iota}{\partial \beta})^{-1}(\frac{\partial u}{\partial \beta } + \iota\frac{\partial v}{\partial \beta } + \frac{\partial \iota}{\partial \beta }v) =
\end{equation}
\begin{equation}
(\frac{\partial u}{\partial \alpha } + \iota\frac{\partial v}{\partial \alpha } + \frac{\partial \iota}{\partial \alpha }v)(\frac{\partial\iota}{\partial \alpha})^{-1} +
(\frac{\partial u}{\partial \beta } + \iota\frac{\partial v}{\partial \beta } + \frac{\partial \iota}{\partial \beta }v)(\frac{\partial\iota}{\partial \beta})^{-1}
\end{equation}
then
\begin{equation}
\iota_{\alpha}^{-1}\iota \frac{\partial v}{\partial \alpha} + \iota_{\beta}^{-1}\iota \frac{\partial v}{\partial \beta}=0
\end{equation}

as the dot product of $\iota_{\alpha}^{-1}\iota$ and $\iota_{\beta}^{-1}\iota$ is zero we conclude that

\begin{equation}
\frac{\partial v}{\partial \alpha} = \frac{\partial v}{\partial \beta} = 0
\end{equation}
and from here

\begin{equation}
\frac{\partial u}{\partial \alpha} = \frac{\partial u}{\partial \beta} = 0
\end{equation}

So $f$ is of Class III. The reciprocal is immediate.
%%%%%%%%%%%%%%%%%%%%%%%%%%%%%%%%%%%%%%%%%%%%%%%%%%%%%%%%%%%%%%%%%%%%%%%%%%%%%%%
%                                SECTION V
%%%%%%%%%%%%%%%%%%%%%%%%%%%%%%%%%%%%%%%%%%%%%%%%%%%%%%%%%%%%%%%%%%%%%%%%%%%%%%%

\section*{V. REGULAR IC FUNCTIONS}
We now look for sufficient and necessary conditions that a CI function has to meet in order to be regular, in other words, that $f$ satisfies:
\begin{eqnarray}
\frac{\partial f}{\partial t} + \iota \frac{\partial f}{\partial r} + (\frac{\partial \bar{p}}{\partial \alpha})^{-1} \frac{\partial f}{\partial \alpha} + (\frac{\partial \bar{p}}{\partial \beta})^{-1} \frac{\partial f}{\partial \beta} = 0
\end{eqnarray}
Because the function is CI, it only has one complex component, so:
\begin{eqnarray}
\frac{\partial u}{\partial \alpha} = \frac{\partial u}{\partial \beta} =
\frac{\partial v}{\partial \alpha} = \frac{\partial v}{\partial \beta} =
0
\end{eqnarray}
Which occurs if and only if
\begin{equation}
\frac{\partial f}{\partial t} + \iota \frac{\partial f}{\partial r} = \frac{2v}{r}
\end{equation}
As before, let us call $\tilde f$ the (unique) complex component of the quaternionic function $f$. It is clear that $\tilde f$ has to satisfy the non-analytic condition:
\begin{equation}
\frac{\partial f}{\partial x} + i \frac{\partial f}{\partial y} = \frac{2v}{y}
\end{equation}
We summarize all of this in the following

\newtheorem{prop7}[prop1]{Proposition}
\begin{prop7}
(Rinehart conditions) Let f be a CI and regular function and $\tilde f$ its complex component, then.
\begin{equation}
\frac{\partial f}{\partial_{l} \bar{p}} = 0 \Longleftrightarrow \frac{\partial \tilde{f}}{\partial \bar{z}} = \frac{2v}{y}
\end{equation}
Let f be a Class III quaternionic function and $\tilde f$ its complex component, then
\begin{equation}
\frac{\partial \tilde{f}}{\partial \bar{z}} = 0 \Longleftrightarrow \frac{\partial f}{\partial_{l}\bar{p}} = \frac{-2v}{r}
\end{equation}
\end{prop7}
\textit{Proof} The first affirmation is immediate. For the second affirmation is just enough to note that Class III functions are in particular Class II.

As a consequence we can provide a functional that turns analytical functions onto complex functions that satisfy (72) and therefore will result into regular functions when transformed in quaternionic functions.
\newtheorem{prop8}[prop1]{Proposition}
\begin{prop8}
Let $z = x +iy$ and $f = u + iv$ be a analytical function, defined on the upper complex plane. Let $f'(z)$ be its total derivative evaluated at a point $z$. Define a functional $L$ as:
\begin{equation}
L(f)(z) = \frac{i}{y}f'(z)-i\frac{v(z)}{y^{2}}
\end{equation}
Then $L(f)$ will satisfy (72)
\end{prop8}

%%%%%%%%%%%%%%%%%%%%%%%%%%%%%%%%%%%%%%%%%%%%%%%%%%%%%%%%%%%%%%%%%%%%%%%%%%%%%%%
%                                SECTION VI
%%%%%%%%%%%%%%%%%%%%%%%%%%%%%%%%%%%%%%%%%%%%%%%%%%%%%%%%%%%%%%%%%%%%%%%%%%%%%%%
\section*{VI. THE IMAGINARY DERIVATIVE}

We now investigate the generation of regular functions taking as seed
Class II functions. We present what we shall call the \textit{imaginary derivative}:

\begin{equation}
\frac{\partial}{\partial_{l} \iota} := (\frac{\partial \iota}{\partial \alpha})^{-1} \frac{\partial f}{\partial \alpha} + (\frac{\partial \iota}{\partial \beta})^{-1} \frac{\partial f}{\partial \beta}
\end{equation}

\newtheorem{prop9}[prop1]{Proposition}
\begin{prop9}

Let $f$ be a (left)-Class II function, then

\begin{equation}
\frac{\partial f}{\partial_{l} \iota} = 2v
\end{equation}

\end{prop9}

\textbf{Proof} This is immediate after we remember a Class II function is Class I, so:

\begin{equation}
\frac{\partial f}{\partial _{l} \bar{p}} = (\frac{\partial f}{\partial t} + \iota \frac{\partial f}{\partial r} ) - \frac{1}{r} \frac{\partial}{\partial_{l} \iota} (f) = \frac{-2v}{r}
\end{equation}

We now use the fact that for Class II functions the Fueter operator (and the imaginary derivative) result in a scalar functions, let $f$ be again a (left)-Class II then

\begin{equation}
\frac{\partial }{\partial_{l} \bar{p}} \frac{\partial f }{\partial_{l} \bar{p}} = \frac{\partial }{\partial_{r} \bar{p}} \frac{\partial f}{\partial_{l} \bar{p}}
\end{equation}

After reordering we obtain

\begin{equation}
\frac{\partial }{\partial_{l} \bar{p}} (\frac{\partial f }{\partial_{l} \bar{p}} - \frac{\partial f }{\partial_{r} \bar{p}} ) = 0
\end{equation}

So the \textit{chiral difference} of a Class II function is regular.

%%%%%%%%%%%%%%%%%%%%%%%%%%%%%%%%%%%%%%%%%%%%%%%%%%%%%%%%%%%%%%%%%%%%%%%%%%%%%%%
%                                SECTION VII
%%%%%%%%%%%%%%%%%%%%%%%%%%%%%%%%%%%%%%%%%%%%%%%%%%%%%%%%%%%%%%%%%%%%%%%%%%%%%%%

\section*{VII. LAURENT SERIES}

Let $f$ be a Class I function. We denote by $\iota(\alpha_{0}, \beta_{0})$ a fixed point in $S^2$.  Let $c$ be a point in this upper complex plane. $c$ can be written as
$c_{1}+c_{2} \iota(\alpha_{0}, \beta_{0})$ where $c_{1},c_{2}$ are real numbers. We consider the annulus
\begin{equation}
U(c,s,S, \alpha_{0}, \beta_{0}) := \{ t,r | s < \sqrt{(t-c_{1})^2 + (r-c_{2})^2} < S \}
\end{equation}

As $f$ is of Class I then its restriction to this annulus is a complex holomorfic function. Lets denote
$f(t,r,\alpha_{0},\beta_{0})$ this restriction. In this annulus we represent this complex function by its Laurent series:
\begin{equation}
f(t,r,\alpha_{0},\beta_{0}) = \sum_{n = -\infty}^{\infty} a_{n}(\alpha_{0},\beta_{0})(t+r\iota(\alpha_{0},\beta_{0})-c_{1}-c_{2} \iota(\alpha_{0}, \beta_{0}))^{n}
\end{equation}
Now let $V(\alpha,\beta)$ denote a connected open set in $S^{2}$. We consider an open set $W$ defined in the quaternionic space as the cartesian product of these two open sets :

\begin{equation}
W := U(c,s,S, \alpha_{0}, \beta_{0}) \rtimes V(\alpha,\beta)
\end{equation}

We observe that the function defined as:

\begin{equation}
\sum_{n = -\infty}^{\infty} a_{n}(\alpha,\beta)(t+r\iota(\alpha,\beta)-c_{1}-c_{2} \iota(\alpha, \beta))^{n}
\end{equation}

concides with $f$ in W. And thus we conclude $f$ can be represented by this Laurent series in the open set $W$. However, if $f$ is of Class II we have the following

\newtheorem{prop19}[prop1]{Proposition}
\begin{prop19}
Let $f$ be a Class II function defined in a open set W as described above. Then is this open set:

\begin{equation}
f(p) = \sum_{n = -\infty}^{\infty} a_{n}(\alpha,\beta)(t+r\iota(\alpha,\beta)-c_{1}-c_{2} \iota(\alpha, \beta))^{n}
\end{equation}
and $a_{n}(\alpha,\beta)$ is of Class II for all n.
\end{prop19}

\textit{Proof} We have already showed the construction of this Laurent series. Lets show first that
$a_{0}(\alpha,\beta)$ is of Class II, lets write $a_{0}(\alpha,\beta) = a_{1,0}(\alpha,\beta)+\iota a_{2,0}(\alpha,\beta)$ where $a_{1,0}(\alpha,\beta)$ and $a_{2,0}(\alpha,\beta)$ are real functions. We write
\begin{equation}
a_{0}(\alpha,\beta) =\frac{1}{2\pi \iota(\alpha,\beta)} \oint_{\gamma} \frac{f(t +r \iota(\alpha,\beta))}{(t+r\iota(\alpha,\beta)-c_{1}-c_{2} \iota(\alpha, \beta))}(dt + \iota(\alpha,\beta)dr)
\end{equation}
Where $\gamma$ is a path that lies in the annulus for each $(\alpha,\beta)$ in W.
We now apply the imaginary derivative to both sides:

\begin{equation}
\frac{\partial}{\partial \iota} (a_{0}(\alpha,\beta)) = \frac{\partial}{\partial \iota} (\frac{1}{2\pi \iota(\alpha,\beta)} \oint_{\gamma} \frac{f(t +r \iota(\alpha,\beta))}{(t+r\iota(\alpha,\beta)-c_{1}-c_{2} \iota(\alpha, \beta))}(dt + \iota(\alpha,\beta)dr))
\end{equation}
As
\begin{equation}
\frac{1}{2\pi \iota(\alpha,\beta)}\frac{f(t +r \iota(\alpha,\beta))}{(t+r\iota(\alpha,\beta)-c_{1}-c_{2} \iota(\alpha, \beta))}
\end{equation}
is a Class II function we conclude that the imaginary derivative of this function is its imaginary part times two. The integral of this imaginary part is exactly the imaginary part of $a_{0}(\alpha,\beta)$ times two. So  $a_{0}(\alpha,\beta)$ is of Class II. We repeat this argument with $f(p)p^{n}$ to prove that $a_{n}(\iota(\alpha,\beta))$ is a Class II function, This is possible because $f(p)p^{n}$ is a Class II function.

We can now provide a functional that will turn a left-Class II function into a right-Class II function defined on a open set W.

\newtheorem{prop21}[prop1]{Proposition}
\begin{prop21}
(The Mirror operator) Let $f$ be a left-Class II function defined in a open set W as previously described, then 
\begin{equation}
M(f)(p) = \bar{f}(\bar{p})
\end{equation}
\end{prop21}
\textit{Proof} Assume $f$ is left regular. We write $f$ with its Laurent series.
\begin{equation}
f(p) = \sum_{n = -\infty}^{\infty} a_{n}(\alpha,\beta)(t+r\iota(\alpha,\beta)-c_{1}-c_{2} \iota(\alpha, \beta))^{n}
\end{equation}
We apply the operator to both sides:
\begin{equation}
\bar{f}(\bar{p}) = \sum_{n = -\infty}^{\infty} \bar{a_{n}}(\alpha,\beta)(t+r\iota(\alpha,\beta)-c_{1}+c_{2} \iota(\alpha, \beta))^{n}
\end{equation}
which is a series centered in $c_{1}-c_{2} \iota(\alpha, \beta)$. This series is a right-Class II function as 
$(t+r\iota(\alpha,\beta)-c_{1}+c_{2} \iota(\alpha, \beta))^{n}$ is central and $\bar{a_{n}}(\alpha,\beta)$ is a right-Class II function.

\section*{IX.ACKNOWLEDGMENTS}
The initial ideas that lead to this article were suggested by A. Sudbery. Carlos Duran helped with critics and suggestions. Nir Cohen pointed out for the author the existence of a center in the Fueter regularity. The author was partially financed by a FAPESP scholarship.
%%%%%%%%%%%%%%%%%%%%%%%%%%%%%%%%%%%%%%%%%%%%%%%%%%%%%%%%%%%%%%%%%%%%%%%%%%%%%%%
%                                REFERENCES
%%%%%%%%%%%%%%%%%%%%%%%%%%%%%%%%%%%%%%%%%%%%%%%%%%%%%%%%%%%%%%%%%%%%%%%%%%%%%%%


\begin{thebibliography}{99}


\bibitem{H43}
W.~R.~Hamilton,
          {\sl Elements of Quaternions} (Chelsea Publishing Co., N.Y., 1969).
\bibitem{Sudbery}
A. Sudbery,
          {\em Quaternionic Analysis}
          Math. Proc. Cambridge Phil. Soc., Vol. 85, pp. 199-225;
          (1979).
\bibitem{F35}
R.~Fueter, Die funktionentheorie der differentialgleichungen $\Delta u =0$ und
$\Delta \Delta u=0$ mit vier variablen,
          Comment.~Math.~Helv.~{\bf 7}, 307-330 (1935);
          {\em ibidem} {\bf 8}, 371 (1936)

\bibitem{Rinehart}
R. F. Rinehart. Elements of a theory of intrinsic functions on algebras. Duke Math. J. 27, no. 1 (1960), 1Ð19

\bibitem{Cul}
C. G. Cullen, An Integral Theorem for Analytic Intrinsic functions on quaternions.
	Duke Math. J. 32, 139-148 (1965).
\bibitem{Deav}
Deavours, C.A., The Quaternion Calculus.
	 The American Mathematical Monthly, Vol. 80, No 9, pp. 995-1008 (1973).
\bibitem{Leut}
Eriksson, Sirkka-Liisa; Leutwiler, Heinz.,
Hypermonogenic functions and their Cauchy-type theorems.  Advances in analysis and geometry, pp 97-112
	Trends Math, Birkhauser, Basel (2004).
\bibitem{deleo}
S. De Leo and P. Rotelli, Quaternion Analiticity. 
Applied Mathematics Letters 16, 1077-1081 (2003)
\end{thebibliography}
\end{document}